\documentclass[lettersize,journal]{IEEEtran}
\usepackage{cite}
\usepackage{amsmath,amssymb,amsfonts}
\usepackage{algorithmic}
\usepackage{graphicx}
\usepackage{algorithm,algorithmic}
\usepackage{hyperref}
\hypersetup{hidelinks=true}
\usepackage{xcolor}
\usepackage{textcomp}
\usepackage{multirow} 
\usepackage{booktabs} 
\usepackage{array}    
\usepackage{graphicx}
\usepackage[caption=false,font=footnotesize,labelfont=rm,textfont=rm]{subfig}

\newtheorem{lemma}{Lemma}
\newtheorem{theorem}{Theorem}
\begin{document}
	
\title{A New Perspective on Optimal Power Flow: From Control Theory to Numerical Optimization}
	
	\author{Yeming Xu, Qi Liu, Liping Zhang,  Huanshui Zhang, and Minyue Fu
	\thanks{This work was supported by the Original Exploratory Program Project of the National Natural Science Foundation of China (62450004), the Joint Funds of the National Natural Science Foundation of China (U23A20325), and the National Natural Science Foundation of China (52307115). (Corresponding author: Huanshui Zhang)}
	
	\thanks{Yeming Xu, Qi Liu, Liping Zhang, and Huanshui Zhang are with the College
		of Electrical Engineering and Automation, Shandong University of Science and Technology, Qingdao 266590, China (e-mail: ymxu2022@163.com, qi\_liu@sdust.edu.cn, lpzhang1020@sdust.edu.cn, hszhang@sdu.edu.cn). Huanshui Zhang is also
		with the School of Control Science and Engineering, Shandong University,
		Jinan, Shandong 250061, China.}
	
	\thanks{Minyue Fu is with the School of Automation, Guangdong University of
	Technology, Guangzhou 510006, China (e-mail: fumy@sustech.edu.cn).}}

	
	
	\maketitle
	
	\begin{abstract}
The optimal power flow (OPF) problem is a classical non-convex optimization problem with nonlinear constraints and multiple potential solutions. Existing methods often fail to find the global optimal solution, motivating the development of a method to identify multiple local solutions. This paper applies the optimal control principle to solve the OPF problem, aiming to systematically explore all feasible solutions within a single run.
We first handle the constraints in the OPF problem using the penalty method and the Lagrange multiplier, and then reformulate it as an optimal control problem. The connection between the two problems is established by the maximum principle. 
Moreover, a numerical algorithm is proposed to solve this optimal control problem. The state trajectory of the system can rapidly and stably converge to the solution of the OPF problem. By designing the initial control sequence, the proposed method can successfully find multiple solutions and escape local solutions. Finally, simulations were conducted on several representative test systems to validate the potential of the proposed method in identifying multiple solutions for non-convex problems.
	\end{abstract}

	\section{Introduction}
The optimal power flow (OPF) problem aims to find a steady state operating point of a power system that minimizes the cost of electric power generation while satisfying operational constraints \cite{carpentier1962contribution,dommel2007optimal}. 
Due to the nonlinear constraint of the power flow equations, the OPF problem is a non-convex optimization problem and may contain many solutions \cite{kumari2010enhanced,bakirtzis2002optimal,bukhsh2013local,hiskens2002exploring,wang2022necessary}.

The solution of the OPF problem began with the gradient method proposed by Dommel and Tinney in 1968 \cite{dommel2007optimal}. 
To effectively handle the constraints, classical approaches such as the Lagrange multiplier method and the penalty function method have been extensively explored \cite{el2007method,qiu2009literature,osman2004solution}.
Building on these foundations, the interior point method (IPM), which introduces slack variables and incorporates a logarithmic barrier function, has become a mainstream approach due to its fast convergence rate \cite{momoh2002review,forsgren2002interior,capitanescu2007interior}.
However, when dealing with highly nonlinear objective functions or constraints, the effectiveness of these traditional methods gradually decreases, showing limitations and being prone to getting stuck in local optimal solutions. \cite{ebeed2018optimal,abido2002optimal,bukhsh2013local,sojoudi2012physics}.

Heuristic methods have been applied to power systems for their ability to handle complex problems \cite{michalewicz1996evolutionary}. Among them, genetic algorithms (GAs) have attracted considerable attention \cite{devaraj2005genetic,younes2007optimal,kumari2010enhanced,lai1997improved,numnonda1999optimal}. However, these methods often suffer from premature convergence in OPF problems \cite{csen2025optimal}, making it difficult to capture all possible solutions. 

In our recent work, we have successfully applied the optimal control principle (OCP) to optimization problems and achieved notable results \cite{zhang2024optimization,wang2025optimization,Lv2025optimization}. Motivated by the characteristic of optimal controllers to drive system states rapidly and stably to equilibrium, we incorporate OCP into finding the solution of non-convex problems. This paper explores the application of the OCP in the OPF problem and proposes a novel optimization method. This method demonstrates outstanding performance in solving highly nonlinear and non-convex problems. Specifically, the main contributions of this work are as follows:
\begin{enumerate}
	\item The proposed method handles inequality constraints using the penalty method and directly solves for the Lagrange multipliers to ensure satisfaction of the equality constraints.	Compared with approaches that rely on power flow calculations to satisfy equality constraints during the optimization of the OPF problem \cite{abaci2016differential,mahdad2015blackout,zhu2023mathematical}, this strategy avoids complex calculations and 
	effectively prevents the objective function from becoming ill-conditioned as in traditional penalty methods, thereby enhancing the numerical stability of the algorithm used to solve this problem.
	\item A solution method for the OPF problem based on the OCP is proposed. The OPF problem is transformed into a discrete-time optimal control problem, where the optimal controller adjusts the first-order difference equations to minimize a cost function composed of control energy and the objective function. It is shown that the steady state of the system in the optimal control problem corresponds to the stationary point of the original problem. Furthermore, an optimization algorithm with superlinear convergence rate is proposed to solve this optimal control problem. Compared with the traditional methods, the proposed method is guaranteed to converge rapidly and stably.
	\item More importantly, it can effectively escape local solutions of the OPF problem and identify multiple OPF solutions within a single execution by appropriately designing the initial control sequence. Numerical experiments further confirm the effectiveness.
\end{enumerate}

	\section{ACOPF Model and Constraint Handling}\label{s2}
	This section introduces the AC OPF formulation in power systems and the associated operational constraints.
	
	\subsection{Problem Formulation}
	The AC OPF problem can be formulated as
	\begin{align}
		\min \quad &F_{obj} = \sum\limits_{g \in \mathcal{N}_G} {{f_g}({P_{G_g}})}, \tag{1a}\label{1a}\\
		\text{s.t.}	\quad &P_{G_g,\min }^{} \le P_{G_g}^{} \le P_{G_g,\max }^{},\tag{1b}\label{1e}\\
		\quad &Q_{G_g,\min }^{} \le Q_{G_g}^{} \le Q_{G_g,\max }^{}, g \in \mathcal{N}_G,\tag{1c}\label{1f}\\
		\quad &{V_{b,\min }} \le {V_b} \le {V_{b,\max }},\tag{1d}\label{1d}\\
		\quad &{P_b} = {P_{G_b}} - {P_{D_b}},\tag{1e}\label{1b}\\
		\quad &{Q_b} ={Q_{G_b}} - {Q_{D_b}},b \in \mathcal{N}_B,\tag{1f}\label{1c}
	\end{align}where $F_{obj}$ is the objective function; $\mathcal{N}_B$ denotes the set of all buses, and $\mathcal{N}_G$ is the set of buses where generators are located. $P_{G_b}$ and $Q_{G_b}$ represent the active and reactive power outputs of the generator at bus~$b$, while $P_{D_b}$ and $Q_{D_b}$ denote the corresponding active and reactive power demands. 
	For buses without generators, i.e., $b \in \mathcal{N}_B \setminus \mathcal{N}_G$, we set $P_{G_b}=Q_{G_b}=0$.
	The voltage magnitude, net active and reactive power injections at bus $b$ are denoted by $V_b$, $P_b$, and $Q_b$, respectively.
	The power flow equations in \eqref{1b} and \eqref{1c} can be expressed as 
	\begin{align}
		P_b &= V_b\sum_{j\in\mathcal{N}_B} V_j\left(G_{bj}\cos\theta_{bj} + B_{bj}\sin\theta_{bj}\right),  \nonumber\\
		Q_b &= V_b\sum_{j\in\mathcal{N}_B} V_j\left(G_{bj}\sin\theta_{bj} - B_{bj}\cos\theta_{bj}\right), \nonumber
	\end{align}where $\theta_{bj} = \theta_b - \theta_j$, and $\theta_b$ denote the voltage angle at bus $b$; $G_{bj}$ and $B_{bj}$ denote the series conductance and susceptance between buses $b$ and $j$, respectively. 
	\setcounter{equation}{1}
		In this paper, the generation cost function is selected as the optimization objective. 


	\subsection{Handling Constraint} 
	For the inequality constraints in problem (1), we employ the penalty method by directly incorporating penalty terms corresponding to \eqref{1e} - \eqref{1d} into the objective function $F_{obj}$. The resulting equality-constrained optimization problem is then expressed as
	\begin{equation}
		\begin{aligned}
			\min \quad & F_{\mathrm{pen}} = F_{\mathrm{obj}} + P_{\mathrm{pen}}, \\
			\text{s.t.} \quad & \text{(1e) and (1f)}, \label{p2}
		\end{aligned}
	\end{equation}where $F_{\mathrm{pen}}$ is the penalized objective function, and the penalty term is formulated as
	\begin{align}
		P_{\mathrm{pen}}
		= &\, \lambda_1 \sum_{b \in \mathcal{N}_B} \phi(V_b, V_{b,\min}, V_{b,\max}) \nonumber\\
		&+ \lambda_2 \sum_{g \in \mathcal{N}_G} \phi(P_{G_g}, P_{G_g,\min}, P_{G_g,\max}) \nonumber\\
		&+ \lambda_3 \sum_{g \in \mathcal{N}_G} \phi(Q_{G_g}, Q_{G_g,\min}, Q_{G_g,\max}),\nonumber
	\end{align}
	where $\lambda_1$, $\lambda_2$, and $\lambda_3$ are the penalty coefficients.
	
	We define the Lagrangian of problem \eqref{p2} as
	\begin{align}
		L(\boldsymbol{y},\boldsymbol{\mu},\boldsymbol{\nu}) =& F_{\rm pen} + \sum_{b \in \mathcal{N}_B}\mu_b(P_{G_b}-P_{D_b}-P_b)\notag\\
		&+ \sum_{b \in \mathcal{N}_B}\nu_b(Q_{G_b}-Q_{D_b}-Q_b),\label{Lag}
	\end{align}
	where 
	$\boldsymbol{\mu} = [\mu_b]_{b \in \mathcal{N}_B}$ and $\boldsymbol{\nu} = [\nu_b]_{b \in \mathcal{N}_B}$ are the Lagrange multipliers associated with power flow equality constraints.
	The optimization  variable $\boldsymbol{y} = [P_G, Q_G, V, \theta]\in\mathbb{R}^m$ includes the generator active power outputs $P_G = [P_{G_g}]_{g \in \mathcal{N}_G}$,  generator reactive power outputs
	$Q_G = [Q_{G_g}]_{g \in \mathcal{N}_G}$, bus voltage magnitudes $V = [V_b]_{b \in \mathcal{N}_B}$, and bus voltage angles 
	$\theta = [\theta_b]_{b \in \mathcal{N}_B\setminus r}$, where $r$ is the index of the reference bus and $m = 2|\mathcal{N}_G| + 2|\mathcal{N}_B| - 1$. 
	
	For clarity in the subsequent discussion, we report the following well-known lemma. A comprehensive overview can be found in Section 11.3 of \cite{luenberger1984linear}.
	\begin{lemma}[First-order necessary conditions]\label{th1}
		Let $\boldsymbol{y}^\star = [P_G^{\star}, Q_G^{\star}, V^{\star}, \theta^{\star}]$ be a local minimum of $F_{\rm pen}$ subject to the constraints \eqref{1b} and \eqref{1c}. 
		Then, there exist unique Lagrange multipliers $\boldsymbol{\mu}^\star = [\mu_b^\star]_{b \in \mathcal{N}_B}$ and $\boldsymbol{\nu}^\star = [\nu_b^\star]_{b \in \mathcal{N}_B}$ such that
		\begin{align}
			\nabla_{\boldsymbol{y}} L(\boldsymbol{y}^\star, \boldsymbol{\mu}^\star, \boldsymbol{\nu}^\star) = 0.\label{Th1}
		\end{align}
	\end{lemma}
	
	In the rest of the paper, we refer to any triple $(\boldsymbol{y}^\star, \boldsymbol{\mu}^\star, \boldsymbol{\nu}^\star)$ that satisfies \eqref{Th1} and the constraints \eqref{1b} and \eqref{1c} as a stationary point of problem \eqref{p2}.

	
	\section{Solution Method for ACOPF via Optimal Control}\label{s3}
	\subsection{Design of the Optimal Control Problem}
	As stated in Lemma~\ref{th1}, the stationary points of problem \eqref{p2} satisfy
	\begin{equation}
		\nabla_x L(x)=0, \label{kkt}
	\end{equation}where \(x=[\boldsymbol{y},\boldsymbol{\mu},\boldsymbol{\nu}]^\top\in\mathbb{R}^n\) and $n=m+2|\mathcal{N}_B|$.
	Nonlinear equations \eqref{kkt} are typically addressed via Newton's method
	\begin{equation}
		x_{k+1} = x_k + \alpha_k d_k,\nonumber
	\end{equation}
	where $x_k=[\boldsymbol{y}_k, \boldsymbol{\mu}_k, \boldsymbol{\nu}_k]^{\top}$ represents the vector of the variable $x$ at the 
	$k$-th iteration, $\alpha_k > 0$ is step size and $d_k$ is the update direction obtained by solving 
	\begin{equation}\label{newton}
		H(x_k) d_k = -\nabla_x L(x_k),
	\end{equation} where $H(x)$ is the Jacobian matrix of $\nabla_x L(x)$.
We designed the following optimal control problem 
\begin{align}
	\min_{u^y} \max_{u^z} \quad & J = \sum_{k=0}^N \left( L(x_k) + \frac{1}{2}(u_k^y)^\top R_y u_k^y - \frac{1}{2}(u_k^z)^\top R_z u_k^z \right)\nonumber\\
	&\qquad+ L(x_{N+1}), \tag{7a} \label{ocp1} \\
	\text{s.t.} \quad & x_{k+1} = x_k + \begin{bmatrix} u_k^y \\ u_k^z \end{bmatrix}, \tag{7b} \label{ocp2}
\end{align}where \eqref{ocp2} is the discrete-time system, $x_k=[y_k,z_k]^\top$, $z_k=[\boldsymbol{\mu}_k,\boldsymbol{\nu}_k]^\top$ is the system state,  $J$ is the cost function, and $u_k= [ (u_k^y)^\top, (u_k^z)^\top ]^\top$ is the control input. The parameter $N$ represents the time horizon, while \(R_y\in\mathbb{S}_{++}^{m}\) and \(R_z\in\mathbb{S}_{++}^{n-m}\) are control weighting matrices. 
The control sequence ${u_k}$ in problem (7) is designed to drive the state trajectory $\{x_k\}$ toward the stationary point of $L(x)$, as described below.\setcounter{equation}{7}
	
	\begin{theorem}\label{le1}
		 The state trajectory $\{x_k\}$ generated by optimal control problem~(7) converges to the stationary point of problem~\eqref{p2}.
	\end{theorem}
	\begin{IEEEproof}
		The maximum principle \cite{naidu2018optimal} is applied in problem (7) to derive the following equations:
		\begin{align}
			p_{k-1} =& p_k + \nabla_x L(x_k),\label{bw}\\
			0 =& R u_k + p_k,\label{eq}
		\end{align}	where~\eqref{bw} is the costate equation with the terminal condition $p_N = \nabla_x L(x_{N+1})$, $R=blkdiag(R_y,-R_z)$, and~\eqref{eq} is the equilibrium condition. 
		It follows from \eqref{ocp2} that as the system approaches steady state, for any given tolerance $\epsilon > 0$, there exists an integer $k_0 > 0$ such that $||u_k|| < \epsilon$ for all $k > k_0$. Combining with~\eqref{eq}, we obtain $||p_k|| < \epsilon ||R||$ for all $k > k_0$.  It then follows from~\eqref{bw} that
		\begin{equation}
			||\nabla_x L(x_k)|| < 2\epsilon ||R||, \quad \forall k > k_0 + 1,\nonumber
		\end{equation}
		which implies that $\{x_k\}$  converges to the stationary point of problem \eqref{p2}. The proof is now completed.
	\end{IEEEproof}


	\subsection{Design of OCP-based Method}\label{OCPp}

	We derive the gradient and Hessian matrix of $J$ with respect to $\{u_k\}$ to clearly elucidate the numerical algorithm.
		The gradient and Hessian matrix of $J$ with respect to $\{u_k\}$ are given by
		\begin{equation}
			\begin{aligned}
				\nabla_U J
				=(I_{N+1} \otimes R)U+ \tilde L\left[ {\begin{array}{*{20}{c}}
						{\nabla_x L({x_1})}\\
						\vdots \\
						{\nabla_x L({x_{N + 1}})}
				\end{array}} \right],\label{gradient}
			\end{aligned}
		\end{equation}
		\begin{equation}
			\begin{aligned}
				\nabla _U^2J =I_{N+1} \otimes R +\tilde LD\tilde L^\top,\label{Hessian}
			\end{aligned}
		\end{equation}where \( U = [u_0^\top, \dots, u_N^\top]^\top \), \(\tilde L = T_{N+1} \otimes I_n \), \( T_{N+1} \in \mathbb{R}^{{(N+1)} \times {(N+1)}} \) is an upper-triangular matrix with all entries on and above the diagonal equal to 1, and \( D \) is a block-diagonal matrix whose \( i \)-th diagonal block is given by \(H(x_i) \), \( i = 1, \dots, N+1 \).

 Denote $U_t=[u_0^{t \top}, \dots, u_N^{t \top}]^\top$ as the vector of the variable $U$ at the $t$-th iteration. The iterative scheme of the proposed numerical algorithm is given as 
	\begin{equation}
		\begin{aligned}
			U_{t + 1} &= U_t  - \alpha_t g_t(U_t),\\
			g_j(U_t) &= \beta_t \left( \nabla_U J(U_t) + M g_{j-1}(U_t) \right), j=1,\dots,t,\label{ocp3}
		\end{aligned}
	\end{equation}
	where $\beta_t=(M+ \nabla_U^2 J(U_t))^{-1} $, $g_0(U_t)=\beta_t \nabla_U J(U_t)$, $M \in\mathbb{S}_+^{n(N+1)}$ is a adjustable matrix.
	The convergence rate analysis of the proposed numerical algorithm is as follows.
	\begin{theorem}
		Assume that  $\nabla_U J(U^*)=0$ and $\nabla^2_U J(U^*)\in\mathbb{S}^{++}_n$. If the sequence $\{x_k\}$ generated by \eqref{ocp3} converges and the step size is chosen as $\alpha_t = 1 - \delta^{t+1}$ with $\delta \in [0,1)$, then it admits superlinear convergence.	
	\end{theorem}
	\begin{IEEEproof}
		By direct derivation, we obtain
		\begin{align}
			g_{t}(U^*)=&\big[I_{n(N+1)}-\big((M+\nabla^2_U J(U^*))^{-1}M\big)^{t+1}\big]\nonumber\\
			&\times\nabla^2_U J(U^*)^{-1}\nabla_U J(U^*)\nonumber\\
			=&\boldsymbol{0}_{n(N+1)}, \nonumber\\
			\nabla_U g_{t}(U^*)=&I_{n(N+1)}-\big((M+\nabla^2_U J(U^*))^{-1}M\big)^{t+1}. \nonumber
		\end{align}
		It follows from \eqref{ocp3} that
		\begin{align}
			&||U_{t+1} - U^*||\notag\\
			=&||U_t-U^*- \alpha_t g_t(U_t)||\notag\\
			=&||U_t-U^*- \alpha_t \big[ g_t(U^*)+\nabla  g_t(U^*)(U_t-U^*)+r_t\big]||\notag\\
			\leq&||U_t-U^*-\alpha_t\nabla g_t(U^*)(U_t-U^*)||+||r_t||\notag\\
			=&||\big[\delta^{t+1}I_{n(N+1)}+\alpha_t\big((M+\nabla^2_U J(U^*))^{-1}M\big)^{t+1}\big](U_t-U^*)||\notag\\
			&+||r_t||\notag\\
			\leq&||\delta^{t+1}I_{n(N+1)}+\alpha_t\big((M+\nabla^2_U J(U^*))^{-1}M\big)^{t+1}||\cdot||U_t - U^*||\notag\\
			&+||r_t||,\notag
		\end{align}where $||r_t||=o(||U_t - U^*||)$ denotes the remainder term.
		Noted that $\rho\big((M+\nabla^2_U J(U^*))^{-1}M\big) < 1$, we obtain $\lim_{t \to \infty} \frac{||U_{t+1} - U^*||}{||U_t - U^*||} = 0$, i.e., \eqref{ocp3} is superlinearly convergent.
		
	\end{IEEEproof}
	
	The numerical algorithm is summarized as Algorithm \ref{alg:optimal control problem}.
	\begin{algorithm}[H]
		\renewcommand{\algorithmicrequire}{\textbf{Input:}}
		\renewcommand{\algorithmicensure}{\textbf{Output:}}
		\caption{OCP-based method.}
		\label{alg:optimal control problem}
		\begin{algorithmic}[1]
			\STATE Initialization: $U_0$, $x_0$, $M\ge0$, $t \leftarrow 0$, $\delta$ and $\varepsilon$		
			\REPEAT
			\STATE Forward update \{$ x_{k}^{t} $\};
			\STATE Calculate $\nabla_U J(U_t)$ and  $\nabla^2_U J(U_t)$ via \eqref{gradient} and \eqref{Hessian};
			\STATE Update $U_{t+1}$ using \eqref{ocp3};
			\STATE $t \leftarrow t + 1$   
			\UNTIL $||\nabla_U J(U_t)|| < \varepsilon $
			\ENSURE  \{$ x_{k}^t $\}		
		\end{algorithmic} 
	\end{algorithm}
	
	Algorithm \ref{alg:optimal control problem} is initialized with an initial control sequence $U_0$.
Consider the case where the initial state is initialized as $x_0=x^*$, with $x^*$ denoting a stationary point of problem \eqref{p2}, and the initial control sequence $U_0$ is chosen as the zero sequence.
In this case, the state trajectory generated by $U_0$ satisfies all the conditions in \eqref{ocp2}, \eqref{bw}, and \eqref{eq}, causing the algorithm to terminate immediately and output $x_k^0=x^*$ for $k=1,\dots,N+1$.
 This means that all states remain fixed at $x^*$. Such a zero-sequence initialization leads to an initial trajectory difficult to explore other potential solutions.
	Therefore, to facilitate a more extensive exploration of the feasible domain, the strategic design of $U_0$ is crucial.
	In the following section, unless otherwise specified, $U_0$ is generated as follows:
	First, the feasible domain of the variable $\boldsymbol{y}$ is partitioned into several subregions, and the central points of these subregions are selected to form a set $\mathcal{S}$, ensuring that  $\mathcal{S}$ contains at least $N+1$ elements.
	Then, for each time step $k = 1,\dots, N+1$, $\boldsymbol{y}_k^0$ is assigned by sequentially taking elements from $\mathcal{S}$, while $\boldsymbol{y}_0^0$ is fixed as the given initial value $\boldsymbol{y}_0$. 
	The Lagrange multipliers $\boldsymbol{\mu}^0_k$ and $\boldsymbol{\nu}^0_k$ are initialized as zero vectors for all index $k$. Finally, $U_0$ can be computed from $u_k^0 = x_{k+1}^0 - x_k^0$.

	\section{Numerical Examples}\label{s4}
The effectiveness of the proposed method is evaluated through numerical simulations on several benchmark multimodal ACOPF systems from \cite{bukhsh2013local}, including WB2, WB3, and WB5. 
	
	\subsection{Simulation Configuration}
	We selected two benchmark solvers, i.e., the MIPS solver based on the IPM and the fmincon solver based on SQP, to compare with our proposed method.
	 The initial values for both solvers are sequentially set as $\boldsymbol{y}_k^0$ for $k = 0, \dots, N+1$.
 $R_y = 10^{-10} I_m$ and $R_z = 10^{-10} I_{n-m}$ are is used in all test systems. To obtain stationary points, the matrix $M$ is constructed as a block-diagonal matrix with $N+1$ identical blocks, each being an $n \times n$ diagonal matrix. In each block, the first $m$ diagonal entries are set to $c$, while the remaining diagonal entries are zero. 
	  We set $\varepsilon = 10^{-8}$ as the tolerance for all algorithms. 
	 
	\subsection{2-, 3-, and 5-bus Systems}\label{case1}
	The test systems WB2, WB3, and WB5 are employed to validate the effectiveness of the proposed method in finding multiple solutions in this subsection.
	For the two-bus system WB2 and the three-bus system WB3, the local and global optimal solutions are very close, with only minor differences in their objective function values. In contrast, the five-bus system WB5 includes two generators with distinct costs, leading to a larger disparity in objective function values. All three systems have two feasible solutions, denoted as $S_1$ and $S_2$. The solutions for each test system are summarized in Table~\ref{tab1}.
	\begin{table}[htbp]
		\centering
		\caption{ACOPF Solutions for WB2, WB3, and WB5 Test Systems}
		\renewcommand{\arraystretch}{1.0}
		\setlength{\tabcolsep}{3pt}
		\begin{tabular}{ccccccc}
			\toprule
			\textbf{WB2} & $F_{\mathrm{obj}}$ & Bus & $V_b$ (p.u.) & $\theta_b$ (deg) & $P_{G_b}$ (MW) & $Q_{G_b}$ (Mvar) \\
			\midrule
			\multirow{2}{*}{$S_1$} & \multirow{2}{*}{877.78} & 1 & 0.952 & 0.00 & 438.89 & 94.44 \\
			&                         & 2 & 1.050 & -57.14 & - & - \\
			\hline
			\multirow{2}{*}{$S_2$} & \multirow{2}{*}{905.73} & 1 & 0.950 & 0.00 & 452.86 & 164.32 \\
			&                         & 2 & 0.976 & -64.94 & - & - \\
			\bottomrule	
		\end{tabular}
		
		\vspace{1mm}		
		\begin{tabular}{ccccccc}
			\toprule
			\textbf{WB3} & $F_{\mathrm{obj}}$ & Bus & $V_b$ (p.u.) & $\theta_b$ (deg) & $P_{G_b}$ (MW) & $Q_{G_b}$ (Mvar) \\
			\midrule
			\multirow{3}{*}{$S_1$} & \multirow{3}{*}{417.25} & 1 & 0.951 & 0.00 & 208.62 & 8.93 \\
			&                         & 2 & 0.950 & -27.25 & - & - \\
			&                         & 3 & 0.981 & -30.36 & - & - \\
			\hline
			\multirow{3}{*}{$S_2$} & \multirow{3}{*}{418.14} & 1 & 0.950 & 0.00 & 209.07 & -20.79 \\
			&                         & 2 & 1.011 & -26.36 & - & - \\
			&                         & 3 & 1.050 & -29.23 & - & - \\
			\bottomrule
		\end{tabular}
		
		\vspace{1mm}
		\begin{tabular}{ccccccc}
			\toprule
			\textbf{WB5}& $F_{\mathrm{obj}}$ & Bus & $V_b$ (p.u.) & $\theta_b$ (deg) & $P_{G_b}$ (MW) & $Q_{G_b}$ (Mvar) \\
			\midrule
			\multirow{5}{*}{$S_1$} & \multirow{5}{*}{946.58} & 1 & 1.050 & 0.00 & 181.43 & 124.09 \\
			&                         & 2 & 0.960 & -3.46 & - & - \\
			&                         & 3 & 0.950 & -3.21 & - & - \\
			&                         & 4 & 0.980 & 37.65 & - & - \\
			&                         & 5 & 1.050 & 45.48 & 220.87 & -30.00 \\
			\hline
			\multirow{5}{*}{$S_2$} & \multirow{5}{*}{1082.33} & 1 & 1.010 & 0.00 & 246.05 & 29.27 \\
			&                         & 2 & 0.950 & -6.52 & - & - \\
			&                         & 3 & 0.950 & -6.57 & - & - \\
			&                         & 4 & 1.010 & 358.66 & - & - \\
			&                         & 5 & 1.050 & 363.11 & 98.15 & -30.00 \\
			\bottomrule
		\end{tabular}
		\label{tab1}
	\end{table}
	
\begin{figure}[htbp]
		\centering
		\subfloat[]{\includegraphics[width=2.5in]{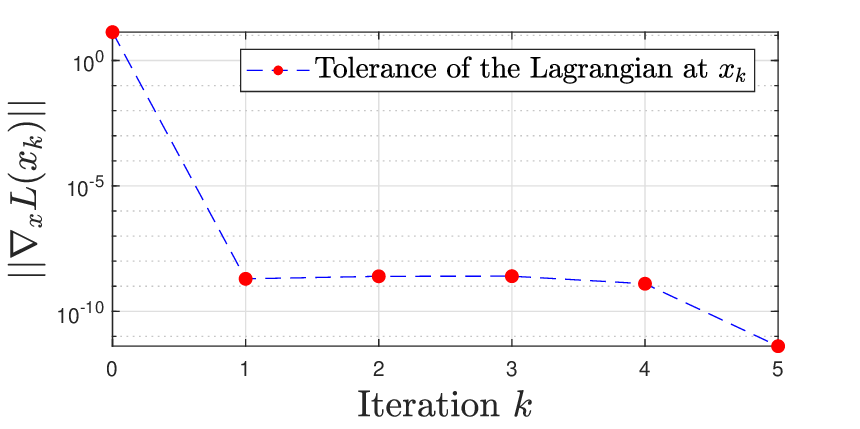}
			\label{WB2_1}}
		\hfill
		\subfloat[]{\includegraphics[width=2.5in]{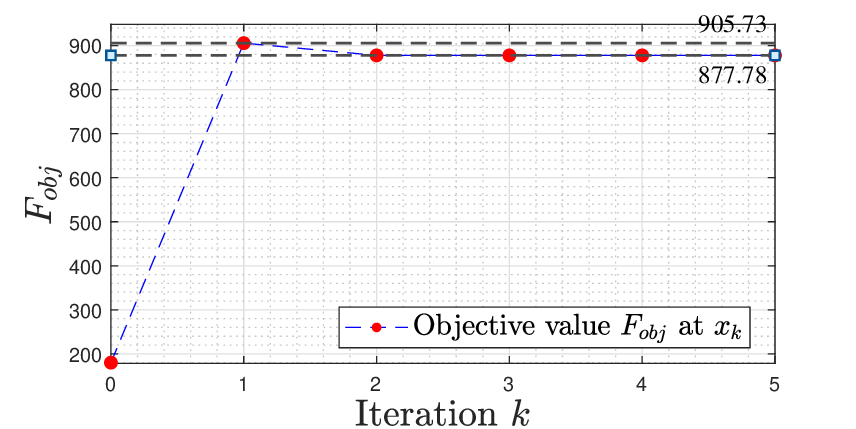}
			\label{WB2_2}}
		\caption{Convergence behavior of the proposed method on WB2. (a) $||\nabla_x L(x)||$ along the state trajectory $\{x_k\}$. (b) Values of $F_{obj}$ along the state trajectory $\{x_k\}$.}
		\label{fig:WB2}
	\end{figure}

	For the case of WB2, we set $N=4$, initial value $\boldsymbol{y}_0=[0.9,-0.4,1,1,2.6]$, $c=0.001$, $\delta=0.9$ and $\lambda_i=10^3,i=1,2,3$ for our method. 
	Fig.~\ref{fig:WB2} shows the results of our proposed method.
	From Fig.~\ref{WB2_1}, it can be seen that $\{x_k\}$ successfully converges, with the corresponding gradient norm of $L(x)$ below $10^{-8}$. Fig.~\ref{WB2_2} illustrates the values of $F_{\rm obj}$ along the state trajectory $\{x_k\}$. It can be observed that $x_2$, $x_3$, $x_4$, and $x_5$ stay at the global optimal solution $S_1$, while $x_1$ stays at the local solution $S_2$. This indicates that all solutions are successfully identified in WB2.
	It should be noted that although $x_0$ corresponds to a smaller $F_{\rm obj}$ for the original problem, it does not satisfy the equality constraints \eqref{1b} and \eqref{1c}. 
	
\begin{figure}[htbp]
		\centering
		\subfloat[]{\includegraphics[width=2.5in]{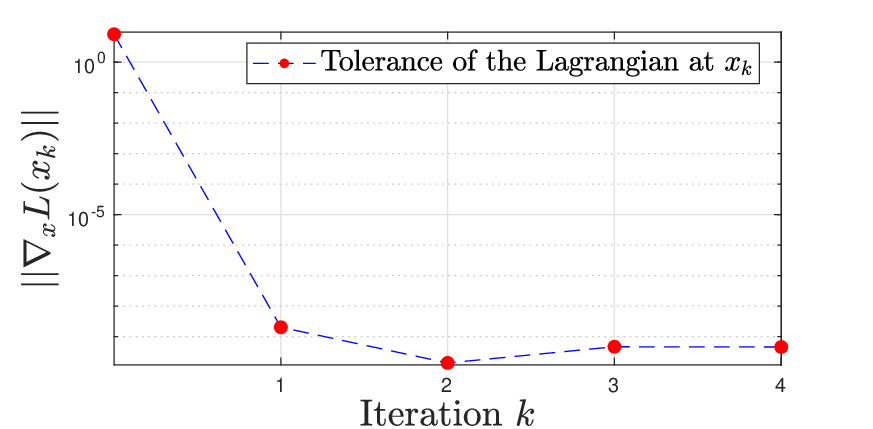}
			\label{WB3_1}}
		\hfill
		\subfloat[]{\includegraphics[width=2.5in]{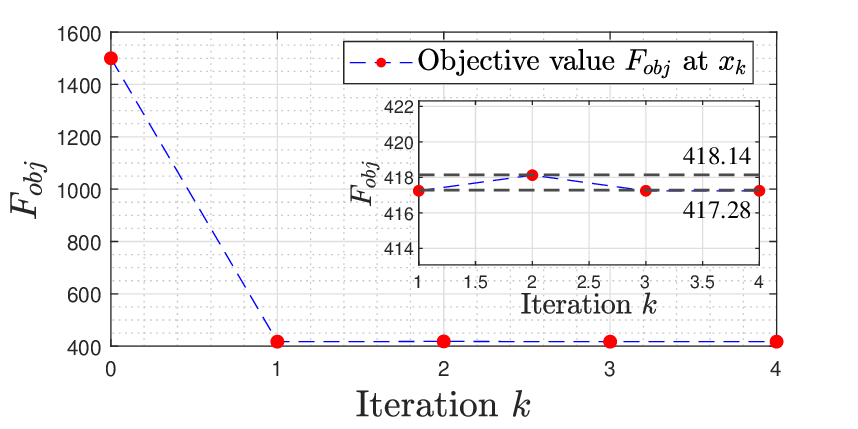}
			\label{WB3_2}}
		\caption{Convergence behavior of the proposed method on WB3. (a) $||\nabla_x L(x)||$ along the state trajectory $\{x_k\}$. (b) Values of $F_{obj}$ along the state trajectory $\{x_k\}$.}
		\label{fig:WB3}
	\end{figure}

	Fig. \ref{fig:WB3} presents the results for the WB3 with $N=3$, $\boldsymbol{y}_0=[7.5,-3,1,1,1,0,0]$, $c=0.0001$, $\delta=0.99$, and $\lambda_i=10^4$ for $i=1,2,3$. As shown in Fig.~\ref{WB3_1}, the state trajectory converges to the stationary points. 
	Fig.~\ref{WB3_2} demonstrates that all solutions are successfully identified, i.e., $x_2$ corresponds to $S_2$, whereas $x_1$, $x_3$, and $x_4$ correspond to $S_1$.
	
\begin{figure}[htbp]
		\centering
		\subfloat[]{\includegraphics[width=2.5in]{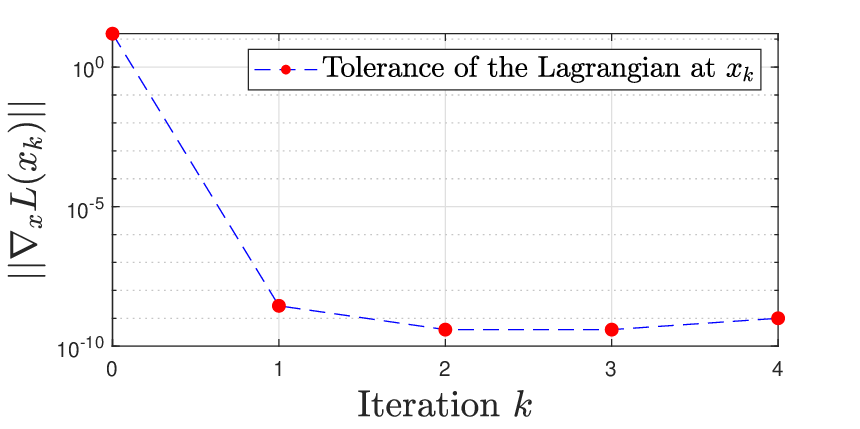}
			\label{WB5_1}}
		\hfill
		\subfloat[]{\includegraphics[width=2.5in]{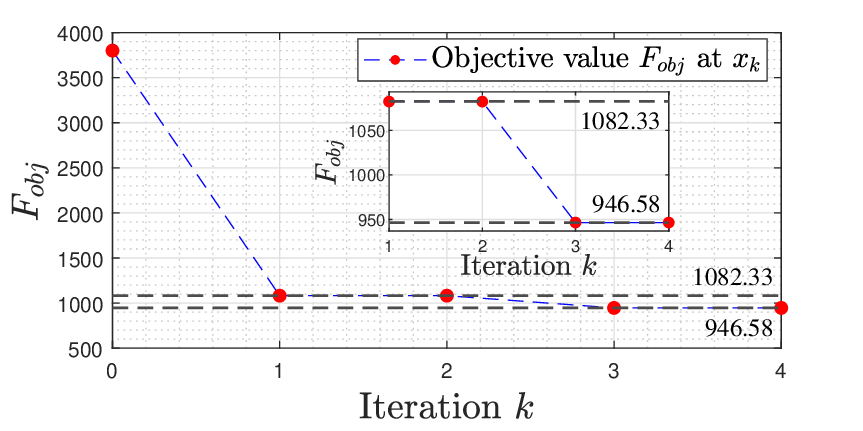}
			\label{WB5_2}}
		\caption{Convergence behavior of the proposed method on WB5. (a) $||\nabla_x L(x)||$ along the state trajectory $\{x_k\}$. (b) Values of $F_{obj}$ along the state trajectory $\{x_k\}$.}
		\label{fig:WB5}
	\end{figure}

	The five-bus system WB5 has two solutions, with the objective value of the local solution being 14\% higher than that of the global optimal solution. 
	We set $N=3$, $\boldsymbol{y}_0=[7.5,8,7.5,8,1,1,1,1,1,0,0,0,0]$, $c=0.01$, $\delta=0.99$, and $\lambda_i=10^3$ for $i=1,2,3$. As we can see from Fig. \ref{WB5_1} and  Fig. \ref{WB5_2}, the state trajectory $\{x_k\}$ generated by the optimal control problem (7) successfully identifies the two solutions of WB5, i.e., $x_3$ and $x_4$ stay at the global optimal solution $S_1$, while $x_1$ and $x_2$ stay at the local solution $S_2$.
	
	As shown in Table~\ref{tab1:total_solutions}, MIPS finds only one solution in WB5 by running with different initial values, and it fails to converge in WB2 and WB3. Meanwhile, fmincon fails to converge for some initial values $\boldsymbol{z}_0$ in WB2, but still identifies two solutions overall. In contrast, it finds only one solution in WB3 and WB5. Our method not only successfully converges stably in all three test systems but also identifies all solutions in a single run. This highlights the effectiveness of the proposed method.
	\begin{table}[htbp]
		\centering
		\caption{Solutions Obtained by Different Solvers in WB2, WB3, and WB5 under Different Initial Values}
		\setlength{\tabcolsep}{6pt}
		\renewcommand{\arraystretch}{1.2}
		\begin{tabular}{c cc cc cc}
			\toprule
			& \multicolumn{2}{c}{WB2} & \multicolumn{2}{c}{WB3} & \multicolumn{2}{c}{WB5} \\
			\cmidrule(lr){2-3} \cmidrule(lr){4-5} \cmidrule(lr){6-7}
		 & MIPS & fmincon & MIPS & fmincon & MIPS & fmincon \\
			\midrule
			$\boldsymbol{y}^0_0$ & N/A & $S_2$ & N/A & $S_1$ & $S_2$ & $S_2$ \\
			$\boldsymbol{y}^0_1$ & N/A & $S_2$ & N/A & $S_1$ & $S_2$ & $S_2$ \\
			$\boldsymbol{y}^0_2$ & N/A & $S_1$ & N/A & $S_1$ & $S_2$ & $S_2$ \\
			$\boldsymbol{y}^0_3$ & N/A & N/A   & N/A & $S_1$ & N/A & $S_2$ \\
			$\boldsymbol{y}^0_4$ & N/A & N/A   & N/A & $S_1$ & $S_2$ & $S_2$ \\
			$\boldsymbol{y}^0_5$ & N/A & $S_1$ & - & -   & -   & -   \\
			\bottomrule
		\end{tabular}		\label{tab1:total_solutions}
		
		\vspace{1mm}
		{\footnotesize \hfill Note: `N/A` indicates the solver failed to converge. \par}
	\end{table}

	\section{Conclusion}\label{s5}
This paper has proposed an optimization method for finding multiple solutions in the OPF problem. Specifically, the original constrained OPF problem has been transformed into an optimal control problem by applying the penalty method, introducing Lagrange multipliers, and leveraging the OCP. The resulting problem has then been solved efficiently using a numerical algorithm with superlinear convergence rate.
To validate the effectiveness of the proposed method, simulations were conducted on several representative test systems. The results have shown that the method successfully finds all solutions and escapes local solutions, demonstrating the potential of the OCP for handling non-convex problems.
	
	\bibliographystyle{IEEEtran}
	\bibliography{reference}

	\vfill
	
\end{document}